\documentclass[letter,11pt]{amsart}
\usepackage{amsfonts,textcomp,amssymb,amsmath,amsthm}
\usepackage{color,ulem}
\usepackage{fullpage}
\usepackage{graphicx}
\usepackage[active]{srcltx}
\usepackage{mathabx}
\usepackage{dsfont,mathrsfs,stmaryrd,wasysym,amsbsy}
\usepackage{enumerate}
\usepackage{IEEEtrantools}
\usepackage{enumitem}   
\usepackage{mathtools}

\newcommand{\inv}[1]{\frac{1}{#1}}
\newcommand{\Xzero}{X_{0-}}
\newcommand{\X}[1]{X_{t}}
\newcommand{\Xminus}[1]{X_{{#1}-}}

\newcommand{\Lminus}[1]{\Lambda_{{#1}-}}

\newcommand{\Lm}[1]{\underline{\Lambda}_{#1}}

\newcommand{\ChiM}{\overline{\chi}}

\newcommand{\PP}{\mathbb{P}}
\newcommand{\EE}{\mathbb{E}}
\newcommand{\Prob}[1]{\mathbb{P}\left(#1\right)}
\newcommand{\E}[1]{\mathbb{E}\left[#1\right]}
\newcommand{\R}{\mathbb{R}}

\newcommand{\f}[2]{\frac{#1}{#2}}

\DeclareMathOperator*{\ess}{ess}

\renewcommand{\inv}[1]{\frac{1}{#1}}

\renewcommand{\epsilon}{\varepsilon}

\newcommand{\sqrttp}{\sqrt{\f{2}{\pi}}}

\newcounter{casenum}
\newcounter{subcasenum}
\addtocounter{subcasenum}{1}

\theoremstyle{plain}
\newtheorem{thm}{Theorem}[section]
\newtheorem{prop}[thm]{Proposition}
\newtheorem{coro}[thm]{Corollary}
\newtheorem{lem}[thm]{Lemma}

\theoremstyle{definition}

\newtheorem{rmk}[thm]{Remark}

\numberwithin{equation}{section}

\begin{document}
	
\title[Well-posedness of the supercooled Stefan problem]{Well-posedness of the supercooled Stefan problem with oscillatory initial conditions}

\author{Scander Mustapha}
\address{Program in Applied \& Computational Mathematics, Princeton University, Princeton, NJ 08544.}
\email{mustapha@princeton.edu}
\author{Mykhaylo Shkolnikov}
\address{ORFE Department, Bendheim Center for Finance, and Program in Applied \& Computational Mathematics, Princeton University, Princeton, NJ 08544.}
\email{mshkolni@gmail.com}
\footnotetext[1]{M.~Shkolnikov is partially supported by the NSF grant DMS-2108680.}

\begin{abstract}
We study the one-phase one-dimensional supercooled Stefan problem with oscillatory initial conditions.~In this context, the global existence of so-called physical solutions has been shown recently in \cite{cuchiero2020propagation}, despite the presence of blow-ups in the freezing rate.~On the other hand, for regular initial conditions, the uniqueness of physical solutions has been established in \cite{delarue2022global}.~Here, we prove the uniqueness of physical solutions for oscillatory initial conditions by a new contraction argument that replaces the local monotonicity condition of \cite{delarue2022global} by an averaging condition.~We verify this weaker condition for fairly general oscillating probability densities, such as the ones given by an almost sure trajectory of $(1 + W_x - \sqrt{2x|\log{|\log{x}|}|})_{+}\wedge 1$ near the origin, where $W$ is a standard Brownian motion.~We also permit typical deterministically constructed oscillating densities, including those of the form $(1 + \sin{1/x})/2$ near the origin.~Finally, we provide an example of oscillating densities for which it is possible to go beyond our main assumption via further complementary arguments.
\end{abstract}

\maketitle

\section{Introduction}

Consider the one-phase one-dimensional supercooled Stefan problem for the heat equation
\begin{equation}\label{PDE}
\begin{cases}
\partial_t u(t,x) = \inv{2}\partial_{xx}u(t,x),\quad x>\Lambda_t,\;\; t>0, \\ 
u(0, x) = f(x),\;\; x\ge 0\quad\textrm{and}\quad u(t,\Lambda_{t}) = 0,\;\; t>0, \\
\dot{\Lambda}_t = \inv{2}\partial_x u(t,x+)|_{x=\Lambda_{t}},\;\; t\ge 0, \\
\Lambda_0=0
\end{cases}
\end{equation}
with a non-negative initial condition $f$.~The unknowns are $u$, the negative of the temperature of a liquid relative to its equilibrium freezing point, as a function of time and space, and the free boundary~$\Lambda_{}$, which encodes the location of a liquid-solid frontier over time.~The temperature is required to solve the heat equation with Dirichlet-type boundary conditions, while the free boundary moves at a speed proportional to the space derivative of the temperature at said boundary (``Stefan condition'').~To ease exposition, we normalize the latent heat coefficient, usually denoted by $\alpha$, to~$1$. 

\medskip

It turns out that, for generic initial conditions, the frontier $\Lambda_{}$ can exhibit jump discontinuities (see, e.g., \cite[Theorem 1.1]{hambly2019mckean}).~A way to circumvent this issue is to restate \eqref{PDE} in a probabilistic form, which allows the definition of global solutions, even in the presence of jump discontinuities.~To wit, let $\Xzero$ be a non-negative random variable with a density $f$, and let $B$ be an independent standard Brownian motion.~The probabilistic reformulation of \eqref{PDE}, first introduced in \cite{delarue2015global, delarue2015particle} for a variant of it, is phrased in terms of the McKean-Vlasov problem
\begin{equation}\label{SDE}
    \begin{cases}
    X_t = \Xzero + B_t - \Lambda_{t},\quad t\ge 0,\\
    \tau := \inf\{t\ge 0:\, X_t\le 0\},\\
    \Lambda_{t} = \Prob{\tau\le t},\quad t\ge 0,
    \end{cases}
\end{equation}
with the unknowns $X = (X_t)_{t\ge 0}$ and $\Lambda_{} = (\Lambda_{t})_{t\ge 0}$.~When $f$ belongs to the Sobolev space $W^1_2([0,\infty))$ and $f(0)=0$, a solution $(X, \Lambda_{})$ of \eqref{SDE} such that  $\dot{\Lambda}\in L^2([0,T])$ for some $T\in(0,\infty)$ gives rise to a solution $u\in W_2^{1, 2}(\{(t,x)\in[0,T]\times[0,\infty)\!:x\ge\Lambda_t\})$ of \eqref{PDE} on $[0,T]$ by taking $u(t, x)\,\mathrm{d}x$ as the law of $(X_t+\Lambda_t)\,\mathbf{1}_{\{\tau>t\}}$ on $(\Lambda_t,\infty)$, for $t\in [0, T]$ (cf.~\cite[proof of Proposition 4.2(b)]{nadtochiy2019particle}).

\medskip

The probabilistic formulation \eqref{SDE} brings out the necessary presence of jump discontinuities in the frontier $\Lambda_{}$ for certain initial data $\Xzero$ (for example, those with $\E{\Xzero} < 1/2$, see \cite[Theorem 1.1]{hambly2019mckean}), as well as the non-uniqueness of the jump sizes $ \Xminus{t} - X_t \coloneqq \lim_{s\uparrow t} X_{s} - X_{t} = \Lambda_{t} - \Lambda_{t-}$ at the instants of discontinuity.~When extending solutions beyond a discontinuity, one must decide how to choose the jump size, which has led to the introduction of the condition
\begin{equation}\label{def:physical}
 \Xminus{t} - X_t = \Lambda_{t} - \Lminus{t} = \inf\big\{x > 0 :\, \Prob{\tau\ge t,\,\Xminus{t}\in (0, x]} < x\big\},\quad t\ge0.
\end{equation}
Solutions of \eqref{SDE} satisfying \eqref{def:physical} are called physical.~It has been shown that \eqref{def:physical} selects the minimal jump sizes a right-continuous solution $\Lambda$ with left limits can have (see \cite[Proposition~1.2]{hambly2019mckean}). The global existence of physical solutions is known under natural assumptions on the initial data, see \cite{cuchiero2020propagation}, where it is proved under the very mild assumption $\E{\Xzero} < \infty$, as well as earlier results in \cite{delarue2015global}, \cite{nadtochiy2019particle}, \cite{nadtochiy2020mean}.~On the other hand, it has been established in~\cite{delarue2022global} that if $\Xzero$ possesses a density $f$ on $[0, \infty)$ that is bounded and changes monotonicity finitely often on compacts intervals, then the physical solution is unique. 

\medskip

This paper develops new arguments that demonstrate uniqueness for oscillatory initial data, which in particular do not fulfill the monotonicity change assumption of \cite{delarue2022global}, though densities fulfilling the latter assumption are also captured by our main theorem.~Oscillatory initial conditions arise frequently when one investigates continuum limits of interacting particle systems. For example, \cite[Remark 1.10]{dembo2019criticality} and \cite[Theorem 1.2]{KoSh} feature initial conditions given by the trajectories of a Brownian motion and a (reparameterized) Brownian bridge, respectively.~We also refer to \cite[Theorem 4]{Lan}, \cite[Theorem 5.4]{FeMe}, \cite[Theorem 4]{CoEi} where the initial conditions even are distributions rather than functions in general.

%For instance, the interacting particle system in the recent work \cite{dembo2019criticality} can be viewed as a space discretized version of \eqref{PDE}, and the initial conditions therein as discretizations of a white noise (and more generally distributional derivatives of H\"older continuous functions) centered around~$1$. Note, however, that on the large time scales studied in \cite{dembo2019criticality} the effective density becomes sufficiently regular to fall under the assumptions of~\cite{delarue2022global} (see \cite[display between (1.23) and (1.24)]{dembo2019criticality}).

\medskip

Consider the question of short time uniqueness for \eqref{SDE}--\eqref{def:physical}, assume that $\Xzero$ has a density~$f$, and let $F$ be the cumulative distribution function (CDF) of $\Xzero$. If $\ess\limsup_{x\downarrow 0} f(x)<1$, there is no jump discontinuity at time $0$ (i.e., $\Lambda_{0} = 0 =:\Lminus{0}$ and $X_0=\Xzero$) for any physical solution, and it is straightforward to prove short time uniqueness.~If $\ess\liminf_{x\downarrow 0} f(x) > 1$, any physical solution must have an initial jump of the size $\Lambda_{0} = \inf\{x > 0\!: F(x) < x\} > 0$, and one can focus on the problem started from $X_{0} = \Xzero - \Lambda_{0}$, with the density $f(x + \Lambda_{0})$.~This new density satisfying necessarily $\ess\liminf_{x\downarrow 0} f(x) \le 1$, we infer that ultimately one needs to investigate the case $\ess\liminf_{x\downarrow 0} f(x) \le 1 \le \ess\limsup_{x\downarrow 0} f(x)$.

\medskip

In \cite[Proposition 5.2]{delarue2022global}, short time uniqueness is shown using a contraction argument, based on the fact that for densities satisfying their monotonicity change assumption, there exists a non-decreasing function $h\!:(0,\infty)\to(0,\infty)$, with $h(0+) = 0$, such that for all $x>0$ sufficiently small:
\begin{equation}\label{condition:mono}
    f(x) \le 1 - h(x).
\end{equation}
The key contribution of this paper is the proof of short time uniqueness for densities oscillating down from $1$, and thus violating \eqref{condition:mono}.~Instead, we introduce an averaging condition:~There exists a non-decreasing function $g\!:(0,\infty)\to(0,\infty)$ such that
\begin{subequations}\label{condition}
\begin{equation}\label{condition:less-one}
 f\le 1,\qquad \int_0^\infty x\,f(x)\,\mathrm{d}x<\infty, \qquad\qquad\qquad\qquad\qquad\qquad\qquad\qquad\qquad\;\;
\end{equation}
% \begin{equation}\label{condition:lower-bounded}
%   \exists x_0>0,\ \exists m > 0,\ \forall x\in [0, x_0],\ f(x)\ge m > 0.
% \end{equation}
\begin{equation}\label{condition:mean}
\exists\,\lambda_0>0\quad\forall\,\lambda\in [0, \lambda_0)\quad \forall\, \mu\in [0, 1]:\quad  \int_{\mu}^{\mu + 1}f(\lambda x)\,\mathrm{d}x \le 1 - g\big(\lambda(\mu+1)\big).
\end{equation}
\end{subequations}
Notice that \eqref{condition:mono} implies 
\begin{equation}
\int_{\mu}^{\mu + 1}f(\lambda x)\,\mathrm{d}x 
\le 1-\int_{\mu}^{\mu + 1} h(\lambda x)\,\mathrm{d}x 
\le 1-\int_{(\mu+1)/2}^{\mu + 1} h(\lambda x)\,\mathrm{d}x 
\le 1-\frac{h(\lambda(\mu+1)/2)}{2},
\end{equation}
i.e., \eqref{condition:mean} with $g(x):=h(x/2)/2$.

\medskip

We are now ready to state our main result. 

\begin{thm}\label{thm:uniqueness}
Let $\Xzero\ge0$ possess a density $f$ that satisfies condition \eqref{condition} with a continuous function $g$.~Then, the physical solution $(X, \Lambda_{})$ of \eqref{SDE} started from $\Xzero$ is unique.
\end{thm}

\begin{rmk}
Our proof of Theorem \ref{thm:uniqueness} (see Section \ref{section:proof-thm}) shows that the solution $(X, \Lambda_{})$ of \eqref{SDE} started from $\Xzero$ is \textit{locally} unique even if one weakens the physicality assumption to $\Lambda_0=0$.
\end{rmk}

In the second part of the article, we provide evidence that condition (\ref{condition}) is natural and non-restrictive, by establishing that it is fulfilled by many oscillating densities, like ones given by sample paths of certain stochastic processes.

\begin{coro}\label{coro:sto}
For almost every fixed sample path of a standard Brownian motion $(W_x)_{x\ge 0}$, the physical solution $(X, \Lambda_{})$ of \eqref{SDE} started from $\Xzero\ge0$ is unique if $\Xzero$ has a density $f$ obeying \eqref{condition:less-one} and such that
\begin{equation}	
f(x) = \big(1 + W_x - \sqrt{2x|\log|\log{x}||}\big)_+\wedge 1,\quad x\in [0,1].
\end{equation}
\end{coro}

\smallskip

We also consider deterministically constructed oscillating densities, including the ones in the next corollary. 

\begin{coro}\label{coro:sinus}
For $\alpha>0$, let $\Xzero\ge0$ be a random variable with the density
 \begin{equation}\label{def:density-sinus-32}
 f(x) = \f{1}{2}\left(1 + \sin\inv{x^{\alpha}}\right),\quad x\in (0, a],
 \end{equation}
 where $a\in (0, \infty)$ is defined by $\int_0^a f(x)\,\mathrm{d}x = 1$.~Then, the physical solution $(X, \Lambda_{})$ of (\ref{SDE}) started from $\Xzero$ is unique.
\end{coro}

\begin{rmk} \label{rmk:oscil}
We note that for the densities $f$ of Corollary \ref{coro:sinus}, the decreasing sequence of solutions to $f(x)=1$ approaches $0$ at an arbitrarily high polynomial rate $n^{-1/\alpha}$. Such oscillatory densities are termed ``pathological'' in \cite[Figure 3.1]{ledger2020uniqueness} due to the difficulty of showing uniqueness for them.
\end{rmk}
% \begin{rmk} 
% \item 1
% \end{enumerate}

\smallskip

The last part the paper exhibits a situation in which it is possible to go beyond condition \eqref{condition} and to establish uniqueness for the supercooled Stefan problem via complementary arguments.
\begin{prop}\label{prop:toy}
Fix a $T\in(0,\infty)$, and let $\Xzero\ge0$ be a random variable with the density
\begin{equation}\label{density:toy}
f(x) = 
\begin{cases}
\alpha_1,\quad x\in \underset{n\ge1}{\bigcup}\,[a_{2n}, a_{2n-1}), \\
\alpha_2,\quad x\in \underset{n\ge1}{\bigcup}\,[a_{2n+1}, a_{2n}),
\end{cases}
\end{equation}
where $0 < \alpha_1<1<\alpha_2$, $a_{2n-1} = r^{n-1}a_1$, $a_{2n} = pr^{n-1}a_1$, and $r=pq$, $p, q\in (0, 1)$.~Then, for any $\alpha_2>1$ close enough to $1$, the physical solution $(X, \Lambda_{})$ of (\ref{SDE}) started from $\Xzero$ is unique on $[0,T]$.
\end{prop}
\begin{rmk}
In contrast to the main theorem (Theorem \ref{thm:uniqueness}), Proposition \ref{prop:toy} is a \textit{local} uniqueness result. In particular, we were unable to verify the monotonicity change assumption of \cite{delarue2022global} at~$T$.  
  \end{rmk}

The rest of the article is structured as follows.~In Section~\ref{section:proof-thm}, we introduce notation and prove Theorem~\ref{thm:uniqueness}.~In Subsection~\ref{subsection:check-sto}, we verify, using functional local laws of the iterated logarithm, that~condition~\eqref{condition} is satisfied by many densities obtained from sample paths of suitable stochastic \\ processes.~In Subsection~\ref{subsection:check-sinus}, we consider oscillating densities constructed from periodic functions. In particular, we deduce Corollaries~\ref{coro:sto},~\ref{coro:sinus} from Theorem~\ref{thm:uniqueness} in Subsections~\ref{subsection:check-sto},~\ref{subsection:check-sinus}, respectively. Finally, Section~\ref{section:toy} is devoted to showing Proposition~\ref{prop:toy}.

\medskip

\noindent\textbf{Acknowledgement.}~We thank Li-Cheng Tsai for bringing the interest in Stefan problems with oscillatory initial conditions to our attention. 

%%%%%%%%%%%%%%%%%%%%%%%%%%
\section{Proof of Theorem~\ref{thm:uniqueness}}\label{section:proof-thm}
%%%%%%%%%%%%%%%%%%%%%%%%%%

Throughout the section, $f$ denotes a density as in Theorem \ref{thm:uniqueness}, and we write $F$ for the associated CDF. We also define the continuous strictly increasing function 
 \begin{equation}
\widetilde{g}:\,[0,\infty)\to[0,\infty),\quad x\mapsto x\,g(x)
 \end{equation}
and set
\begin{equation}
\psi(\lambda, \mu) = \int_{\mu}^{\mu + 1} f(\lambda x)\,\mathrm{d}x, \quad \lambda,\mu\ge 0.
\end{equation}
  
Let $(X,\Lambda)$ be an arbitrary physical solution of \eqref{SDE}.~By \cite[Proposition 2.3]{cuchiero2020propagation}, there exists a minimal solution $(\underline{X},\underline{\Lambda})$ of \eqref{SDE}, namely the unique solution of \eqref{SDE} satisfying
\begin{equation}
\Lm{t} \le \widetilde{\Lambda}_t,\quad t\ge 0,
\end{equation}
for any solution $(\widetilde{X},\widetilde{\Lambda})$ of \eqref{SDE}.~The physicality of $(\underline{X},\underline{\Lambda})$ is ensured by~\cite[Theorem 6.5]{cuchiero2020propagation}.~We further introduce $(Y_t)_{t\ge 0}$, $(Z_t)_{t\ge 0}$ given  respectively by
\begin{eqnarray}
&& Y_t=\sup_{0\le s\le t}(-B_{s} + \Lm{s}),\quad t\ge 0, \\
&& Z_t=\sup_{0\le s\le t}(-B_{s} + \Lambda_{s}),\quad t\ge 0.
\end{eqnarray}
In these terms, the frontiers solve
\begin{eqnarray}
&& \Lm{t} = \PP\big(\inf_{0\le s\le t}(\Xzero + B_s - \Lm{s})\le 0\big) \label{F1}
= \E{F(Y_{t})}, \quad t\ge 0, \\
&& \Lambda_{t} = \E{F(Z_{t})},\quad t\ge 0. \label{F2}
\end{eqnarray}

\smallskip

Our starting point is the following continuous upper bound on the frontier $\Lambda$.

%\begin{lem}\label{lem:prelim-lower-bound}
%There exists a continuous function $(\Chim_{t})_{t\ge 0}$, positive at times $t>0$, such that 
%\begin{equation}
%\Lambda_{t}\ge\Chim_t,\quad t\ge0.
%\end{equation}
%\end{lem}

%\noindent\textbf{Proof.} Thanks to the lower bound $Z_t\ge \sup_{0\le s\le t}(-B_s)$ we have
%\begin{equation}
%\Lambda_{t} = \E{F(Z_t)}\ge \EE\big[F\big(\sup_{0\le s\le t}(-B_s)\big)\big],\quad t\ge 0.
%\end{equation}
%The random variable $\sup_{0\le s\le t}(-B_s)$ equals in distribution to $\sqrt{t}\,|\mathcal{N}|$, where $\mathcal{N}$ is a standard normal random variable.~Therefore,
%\begin{equation}
%\Lambda_{t} \ge \EE[F(\sqrt{t}\,|\mathcal{N}|)]=:\Chim_t,\quad t\ge0.
%\end{equation}
%We check easily that $\Chim$ is continuous. \qed
  
\begin{lem}\label{lem:prelim-upper-bound}
There exist a $T>0$ and a strictly increasing continuous function $(\ChiM_{t})_{t\ge 0}$, with $\ChiM_{0}=0$, such that
\begin{equation}\label{lem:upper-bound}
\Lambda_{t}\le \ChiM_t,\quad t\in [0, T]. 
\end{equation}
\end{lem}

\noindent\textbf{Proof.} For $t\ge 0$, we estimate 
    \begin{equation}
      \begin{split}
        \Lambda_{t} &= \PP\big(\Xzero\le \sup_{0\le s\le t} (-B_s + \Lambda_{s})\big) \\
              &\le \Prob{\Xzero\le \Lambda_{t}} + \PP\Big(\{\Lambda_{t}<\Xzero\} \cap\big\{\Xzero\le \sup_{0\le s\le t} (-B_s + \Lambda_{s})\big\}\Big).
      \end{split}
    \end{equation}
In view of the upper bound $\sup_{0\le s\le t}(-B_s+\Lambda_{s})\le \sup_{0\le s\le t}(-B_s) +\Lambda_{t}$, we find for all $t\ge 0$ that
    \begin{equation}\label{lem:prelim-upper-bound-calc}
      \begin{split}
        \Lambda_{t} - F(\Lambda_{t})&\le \PP\big(\{\Lambda_{t}<\Xzero\} \cap \{\Xzero - \Lambda_{t}\le \sqrt{t}\,|\mathcal{N}|\}\big) \\
        &= \int_{0}^{\infty}\PP(x\le \sqrt{t}\,|\mathcal{N}|)\,\PP(\Xzero- \Lambda_{t} \in \mathrm{d}x) \\ 
        &\le \sqrt{2t/\pi},
      \end{split}
    \end{equation}
where $\mathcal{N}$ is a standard normal random variable and we have used $f\le 1$ and $\E{|\mathcal{N}|}=\sqrt{2/\pi}$.

\medskip

To conclude we apply \eqref{condition:mean} to obtain
    \begin{equation}\label{lem:prelim-upper-bound-calc'}
      \f{F(\lambda)}{\lambda} = \int_0^1f(\lambda x)\,\mathrm{d}x \le 1 - g(\lambda),\quad \lambda \in (0, \lambda_0).
    \end{equation}
Since $\Lambda_{t}$ is right-continuous with $\Lambda_{0} = 0$, there exists a $T>0$ such that $\Lambda_{t} < \lambda_{0}$, $t\in[0,T]$. Putting this together with \eqref{lem:prelim-upper-bound-calc'} and
\eqref{lem:prelim-upper-bound-calc} we get
\begin{equation}
\widetilde{g}(\Lambda_{t}) \le \Lambda_{t} - F(\Lambda_{t}) \le \E{|\mathcal{N}|}\sqrt{t},\quad t\in [0, T].
\end{equation}
The proof is completed by inferring 
    \begin{equation}
      \Lambda_{t}\le \widetilde{g}^{-1}(\E{|\mathcal{N}|}\sqrt{t})=:\ChiM_t,\quad t\in [0, T],
    \end{equation}
with the continuous $\widetilde{g}^{-1}$ satisfying $\widetilde{g}^{-1}(0)=0$. \qed

\medskip

We also need the next lemma.

\begin{lem}\label{lem:prelim-continuity-phi}
Let $(\nu_t)_{t\ge 0}$ be a strictly increasing continuous function, with $\nu_0=0$.~Then, there exists a positive function $\varphi(t,b)$ of $t>0$ and $b>0$ so that
\begin{equation}
P(t,b)\coloneqq\PP\big(\sup_{0\le s\le t}(-B_s+\nu_s)\le b\big) \ge \varphi(t,b).
    \end{equation}
  \end{lem}

\noindent\textbf{Proof.} Fix $t>0$ and $b>0$. If $\nu_t\le b/2$, then
\begin{equation}
P(t, b)\ge \PP\big(\sup_{0\le s\le t}(-B_s)\le b - \nu_t\big).
\end{equation}
Otherwise, $\nu_t>b/2$ and $\tau\coloneqq\tau(b)\coloneqq \nu^{-1}(b/2)<t$. Moreover, for any $\tau'\in (0, \tau]$,
\begin{equation}
\begin{split}
P(t, b)&= \PP\Big(\big\{\sup_{0\le s\le \tau'}(-B_s+\nu_s)\le b\big\} \cap 
\big\{\sup_{\tau'< s\le t}(-B_s+\nu_s)\le b\big\}\Big) \\
&\ge \PP\bigg(\bigg\{\sup_{0\le s\le \tau'}(-B_s)\le \f{b}{2}\bigg\} 
\cap \big\{\sup_{\tau'< s\le t}(-B_s)\le b - \nu_t\big\}\bigg).
	\end{split}
\end{equation}
We conclude by setting
\begin{equation}
\varphi(t, b) 
= \PP\bigg(\bigg\{\sup_{0\le s\le \tau(b)\wedge (\nu_t-b/2)_{+}}(-B_s)\le \f{b}{2}\bigg\} \cap \big\{\sup_{\tau(b)\wedge (\nu_t-b/2)_{+}< s\le t}(-B_s)\le b-\nu_t\big\}\bigg).
\end{equation}
Clearly, $\varphi$ is positive on $(0,\infty)^2$. \qed

\medskip

The following proposition is the key ingredient in our proof of Theorem \ref{thm:uniqueness}.

\begin{prop}\label{prop:iteration}
There exists a function $\Phi\!:(0,T]\times(0,\lambda_0)\to(0,\infty)$ such that %for all continuous functions $h\!:[0,T]\to[0,1]$ with $h_t>0$, $t\in(0, T]$, it holds
    \begin{equation}\label{eq prop:iteration}
      \E{F(Y_t + \lambda) - F(Y_t)} \le (1 - \Phi(t,\lambda))\,\lambda,\quad (t,\lambda)\in(0,T]\times(0,\lambda_0).
    \end{equation}
\end{prop}

\noindent\textbf{Proof.} Let $(t,\lambda)\in(0,T]\times(0,\lambda_0)$.~Then,
%    \begin{equation}
%        \begin{split}
%          \E{F(Y_t + h_t\Lambda_{t}) - F(Y_t)}&= \E{1_{Y_t\le a}\left[F(Y_t + h_t\Lambda_{t}) - F(Y_t)\right]}\\ &+\E{1_{a < Y_t}\left[F(Y_t + h_t\Lambda_{t}) - F(Y_t)\right]}.
%        \end{split}
%    \end{equation}
%\smallskip
%    $F$ being constant equal to one over $[a, \infty]$, the second term of the right-hand side is zero.
%    Moreover, for $t>0$
    \begin{equation}\label{prop:iter-psi-computation}
      \begin{split}
       %\E{1_{Y_t\le a}(F(Y_t+h_t\Lambda_{t}) - F(Y_t)} \le 
       \E{F(Y_t+\lambda) - F(Y_t)}
        &= \E{F\left(\lambda\left(\f{Y_t}{\lambda}+1\right)\!\right) - F\left(\lambda\,\f{Y_t}{\lambda}\right)}\\
        &= \E{\psi\left(\lambda, \f{Y_t}{\lambda}\right)}\lambda.
      \end{split}
    \end{equation}
%Using condition (\ref{condition}), there exists $\lambda_0 > 0$ such that for all $\lambda \in [0, \lambda_0)$ and $\mu\in [0, 1]$
%    \begin{equation}
%      \psi(\lambda, \mu) \le 1 - g(\lambda\mu).
%    \end{equation}
Since $\lambda\in(0,\lambda_0)$, condition \eqref{condition:mean} yields 
 \begin{equation}
      \mathbf{1}_{\{Y_t\le\lambda\}}\,\psi\left(\lambda, \f{Y_t}{\lambda}\right)
       \le \mathbf{1}_{\{Y_t\le\lambda\}}\,(1 - g(Y_t+\lambda)).
    \end{equation}

\smallskip

In view of $\psi \le \sup_{x\ge0} f(x) \le 1 $,
    \begin{equation}\label{where to insert}
      \begin{split}
        \frac{\E{F(Y_t+\lambda) - F(Y_t)}}{\lambda} &\le \E{\mathbf{1}_{\{Y_t\le \lambda\}}\,\psi\left(\lambda, \f{Y_t}{\lambda}\right)}
        + \Prob{Y_t>\lambda} \\
        &\le \E{\mathbf{1}_{\{Y_t\le \lambda\}}(1 - g(Y_t+\lambda))} + \Prob{Y_t>\lambda} \\
        &= 1 - \E{\mathbf{1}_{\{Y_t\le\lambda\}}\,g(Y_t+\lambda)}.
      \end{split}
    \end{equation}

\smallskip

Finally, we use $\sup_{0\le s\le t}(-B_s)\le Y_t = \sup_{0\le s\le t}(-B_s+\Lm{s})$ and $\Lm{t}\le \ChiM_t$, $t\in [0, T]$:
    \begin{equation}
      \begin{split}
        \E{\mathbf{1}_{\{Y_t\le \lambda\}}\,g(Y_t+\lambda)}
        %&\ge \EE\big[\mathbf{1}_{\{Y_t\le h_t\Chim_t\}}\,g(Y_t+h_t\Lambda_{t})\big] \\
        &\ge g(\lambda)\,\PP(Y_t\le\lambda) \\
        &\ge g(\lambda)\,\PP\big(\sup_{0\le s\le t}(-B_s+\ChiM_{s})\le\lambda\big).
      \end{split}
    \end{equation}
Thus, thanks to Lemma~\ref{lem:prelim-continuity-phi},
    \begin{equation}
       \E{\mathbf{1}_{\{Y_t\le\lambda\}}\,g(Y_t+\lambda)}\ge  g(\lambda)
       \,\varphi(t,\lambda)=:\Phi(t,\lambda).
    \end{equation}
Inserting this into \eqref{where to insert} we obtain \eqref{eq prop:iteration}. \qed

\medskip  

We are now ready for the proof of Theorem~\ref{thm:uniqueness}.

\medskip

\noindent\textbf{Proof of Theorem~\ref{thm:uniqueness}.} To start, we fix a $\lambda\in(0,\lambda_0)$ and decrease $T>0$ to ensure $\Lambda_{t}-\Lm{t}\le\lambda$, $t\in[0,T]$, relying on right-continuity. For a $\lambda'\in(0,\lambda]$, suppose $\{t\in[0,T]\!:\Lambda_{t}-\Lm{t}\ge\lambda'\}\neq\emptyset$ and consider $t_{\lambda'}\coloneqq\inf\{t\in[0,T]\!:\Lambda_{t}-\Lm{t}\ge\lambda'\}$. Then, $t_{\lambda'}>0$ by right-continuity, and
\begin{equation}\label{Delta bounds}
0<\lambda' \le \Lambda_{t_{\lambda'}}-\Lm{t_{\lambda'}} = \sup_{0\le t\le t_{\lambda'}} (\Lambda_{t}-\Lm{t})\le\lambda<\lambda_0.
\end{equation}
Therefore, we have
\begin{equation}\label{thm:proof-start}
Z_{t_{\lambda'}} = \sup_{0\le s\le t_{\lambda'}} (-B_s+\Lambda_{s}) 
\le Y_{t_{\lambda'}} + \Lambda_{t_{\lambda'}}-\Lm{t_{\lambda'}}.
\end{equation}
Thus, combining \eqref{F1}, \eqref{F2} and Proposition \ref{prop:iteration} we infer
\begin{equation}
\begin{split}
\Lambda_{t_{\lambda'}}-\Lm{t_{\lambda'}} &= \EE[F(Z_{t_{\lambda'}})-F(Y_{t_{\lambda'}})] \\
&\le \EE[F(Y_{t_{\lambda'}}+\Lambda_{t_{\lambda'}}-\Lm{t_{\lambda'}})-F(Y_{t_{\lambda'}})] \\
&\le (1-\Phi(t_{\lambda'},\Lambda_{t_{\lambda'}}-\Lm{t_{\lambda'}}))\cdot(\Lambda_{t_{\lambda'}}-\Lm{t_{\lambda'}}),
\end{split}
\end{equation}
where $\Phi(t_{\lambda'},\Lambda_{t_{\lambda'}}-\Lm{t_{\lambda'}})>0$.~Hence, $\Lambda_{t_{\lambda'}}-\Lm{t_{\lambda'}}=0$, contradicting \eqref{Delta bounds}.~We readily conclude that $\{t\in[0,T]\!:\Lambda_{t}-\Lm{t}\ge\lambda'\}=\emptyset$, and since $\lambda'\in(0,\lambda]$ was arbitrary, $\Lambda_{t}\le\Lm{t}$, $t\in[0,T]$.~Due to the minimality of $\Lm{}$, it must hold $\Lambda_{t}=\Lm{t}$, $t\in[0,T]$.

%\smallskip
%
%Define $M$ the mapping
%    \begin{equation}
%      M : h\in \mathcal{C}\left([0, T], [0, 1]\right)\to (1 - \phi(\cdot, h))h \in \mathcal{C}\left([0, T], [0, 1]\right).
%    \end{equation}
%Start with $h = 1$ and iterate (\ref{thm:proof-iteration}). Setting $h_n = M^{(n)}(1)$, $n\ge 0$ with have
%    \begin{equation}
%      \begin{cases}
%        (1 - h_n)\Lambda_{t}\le \Lm{t},\ t\in [0, T]\\
%       h_{n+1} \le h_{n}.
%      \end{cases}
%    \end{equation}
%    According to Dini's theorem, there exists $h^{*}\in \mathcal{C}\left([0, t], [0, 1]\right)$ such that
%    \begin{equation}
%      \lim_{n\to \infty}\sup_{0\le s\le T}|h_{n+1}(t) - h^{*}(t)| = 0.
%    \end{equation}
%    Moreover for all $t\in [0, T]$
%    \begin{equation}
%      (1 - \phi(t, h^{*}_t))h^{*}_t = h^{*}_{t}.
%    \end{equation}
%    If for some $t\in (0, T]$, $h^{*}_{t}>0$ then $\phi(t, h^{*}_t) = 0$ which is impossible. In conclusion $h^{*}= 0$ and 
%    \begin{equation}
%      \Lambda_{t}\le \Lm{t},\ t\in [0, T].
%    \end{equation}
%    The minimiality of $\Lm{}$ implies that $\Lambda_{t} = \Lm{t}$, for $t\in [0, T]$.

\medskip

To derive global uniqueness, we let
\begin{equation}
T'\coloneq\inf\{t \ge T:\,\Lambda_{t} \neq \Lm{t}\} \in [T, \infty]
\end{equation}
and suppose that $T'<\infty$. By the definition of $T'$,
\begin{equation}
    X_{T'-} = X_{0-}+B_{T'}-\Lambda_{T'-} =  X_{0-}+B_{T'}-\underline{\Lambda}_{T'-} = \underline{X}_{T'-}
\end{equation}
and $\mathbf{1}_{\{\tau\ge T'\}}=\mathbf{1}_{\{\underline{\tau}\ge T'\}}$, so that $\mathbf{1}_{\{\tau\ge T'\}}\,X_{T'-} = \mathbf{1}_{\{\underline{\tau}\ge T'\}}\,\underline{X}_{T'-}$. Moreover, for all $0 < a < b$,
\begin{equation}
\PP\big(\mathbf{1}_{\{\tau\ge T'\}}\,X_{T'-}\in [a, b]\big) 
= \PP\big(\tau\ge T',\,X_{T'-}\in [a, b]\big)
\le \PP(X_{T'-}\in [a, b]).
\end{equation}
Thus, the right essential limit superior of the density of $\mathbf{1}_{\{\tau\ge T'\}}\,X_{T'-} = \mathbf{1}_{\{\underline{\tau}\ge T'\}}\,\underline{X}_{T'-} $ at $0$ is at most that of $X_{T'-}=\underline{X}_{T'-}$, namely $\EE[f(-B_{T'}+\Lambda_{T'-})]=\EE[f(-B_{T'}+\underline{\Lambda}_{T'-})]$. Since $f\le1$, and $f\equiv0$ on $(-\infty,0)$,
\begin{equation}
\EE[f(-B_{T'}+\Lambda_{T'-})]=\EE[f(-B_{T'}+\underline{\Lambda}_{T'-})]<1. 	
\end{equation}
Consequently, the condition \eqref{condition:mono} is satisfied at $T'-$ and we get $\Lambda\equiv\underline{\Lambda}$ on a non-trivial interval $[T',T'+s]$ by repeating \cite[proof of Proposition~5.2]{delarue2022global}. (Note that the condition \eqref{condition:mono} permits us to apply \cite[Lemma~5.1]{delarue2022global}.) This is the desired contradiction. \qed

%%%%%%%%%%%%%%%%%%%%%%%%
\section{Analysis of specific oscillatory initial conditions} \label{section:check}
%%%%%%%%%%%%%%%%%%%%%%%%
  
%%%%%%%%%%%%%%%%%%%%%%%%  
\subsection{Initial conditions constructed from stochastic processes} \label{subsection:check-sto}
%%%%%%%%%%%%%%%%%%%%%%%%

In the present subsection we illustrate Theorem \ref{thm:uniqueness} on initial conditions obtained from sample paths of stochastic processes. Concretely, we consider initial densities
 \begin{equation}\label{density:sto}
 f(x) = 
 \begin{cases}
 (1 + S_x - \kappa_x)_+\wedge 1,\quad x\in [0, 1], \\
  f_{0}(x),\quad x>1,
 \end{cases}
 \end{equation}
where $(S_x)_{x\ge 0}$ is a stochastic process starting at zero, $(\kappa_x)_{x\ge0}$ is a function with $\kappa_0=0$, and the (random) extension $f_0\!:(1,\infty)\to[0,1]$ ensures that $\int_0^\infty f(x)\,\mathrm{d}x=1$ and $\int_0^\infty x\,f(x)\,\mathrm{d}x<\infty$. 

\medskip

Our interest lies in processes $S$ and functions $\kappa$ such that, almost surely, $f$ violates the local monotonicity condition~\eqref{condition:mono} but satisfies condition \eqref{condition}.~Clearly, condition \eqref{condition:mono} is violated if $S_x \ge \kappa_x$ for a sequence of $x$'s converging to $0$, that is, 
\begin{equation}\label{limsupgeone}
\limsup_{x\downarrow 0} \f{S_x}{\kappa_x} \ge 1.
\end{equation}
As a guiding example, take $S$ to be a standard Brownian motion and $\kappa_x = \sqrt{2x|\log{|\log{x}|}|}$. Due to Chung's law of the iterated logarithm (LIL), the resulting $f$ violates condition \eqref{condition:mono} almost surely.~On the other hand, using the local Strassen's LIL of \cite{strassen1964invariance}, \cite{gantert1993inversion} we prove below that condition \eqref{condition} is satisfied almost surely.~This result extends to other centered continuous Gaussian processes admitting a local functional LIL as follows.
  
\medskip

Let $(S_x)_{x\in[0,1]}$ be a centered continuous Gaussian process with $S_0=0$ and a covariance function $\Gamma(x,y)=\E{S_xS_y}$ continuous on $[0,1]^2$ and non-degenerate on $(0,1]^2$.~We write $H(\Gamma)$ for the reproducing kernel Hilbert space associated with $\Gamma$. Recall that $H(\Gamma)$ is defined as the completion of the space of finite linear combinations of $\{\Gamma(x,\cdot)\}_{x\in[0, 1]}$ under the norm induced by the inner product $\langle \Gamma(x,\cdot), \Gamma(y, \cdot)\rangle:= \Gamma(x,y)$. Elements $\phi\in H(\Gamma)$ obey $\phi(x) = \langle \phi, \Gamma(x, \cdot)\rangle$ and are continuous functions. Therefore, $H(\Gamma)$ is a subset of the Banach space $C([0,1])$. Moreover, the unit ball 
\begin{equation}
K = \{\phi\in H(\Gamma):\,\langle\phi,\phi\rangle\le 1\}
\end{equation}
is compact in $C([0,1])$ (see, e.g., \cite[Lemma 3]{oodaira1972strassen}).~For technical reasons, we assume throughout that the process $S$ has the scaling property
\begin{equation}\label{scaling-prop}
(S_{r x})_{x\in[0,1]} \stackrel{d}{=} (\sqrt{r}^{\,\alpha_2}S_x)_{x\in[0,1]},\quad r\in(0,1],
\end{equation}
for some $\alpha_2>0$. Under \eqref{scaling-prop}, there exists an $\alpha_{1}>0$ such that
\begin{equation}
\gamma(x)\coloneqq \Gamma(x,x) = \alpha_{1} x^{\alpha_{2}},\quad x\in[0,1].
\end{equation}
For simplicity, we take $\alpha_1=1$, i.e., $\E{S_1^2}=1$.

\medskip

We say that $S$ satisfies a local functional LIL if the following assertion holds for a $\beta\in(0,\infty)$.

\medskip

\noindent\textbf{Assertion.} Almost surely, the set
\begin{equation}\label{def:xi}
\{(\xi^r_x)_{x\in[0,1]}\}_{r\in(0,1]}:= \bigg\{\bigg(\f{S_{rx}}{\beta\sqrt{\gamma(r)\, |\log{|\log{r}||}}}\bigg)_{\!x\in[0,1]}\bigg\}_{r\in(0,1]}
\end{equation}
is relatively compact in $C([0, 1])$, and the set of its limit points as $r\downarrow0$ is given by $K$. In particular, for every continuous functional $\mathcal{I}\!:C([0,1])\to\R$ we have
\begin{equation}
\limsup_{r\downarrow 0}\,\mathcal{I}(\xi^r) = \sup_{\phi\in K}\,\mathcal{I}(\phi)\quad\text{almost surely.}
    \end{equation}

\smallskip

A local functional LIL has been established in the case of a fractional Brownian motion with Hurst exponent $H\in(0,1)$, for which $\Gamma(x,y) = \inv{2}(x^{2H} + y^{2H} - |y - x|^{2H})$ and $\gamma(x) = x^{2H}$ (see~\cite[Example 4.35]{malyarenko2012invariant}).~Further, by taking $\mathcal{I}(\phi) = \phi_1$ one derives the usual LIL, so that the density $f$ in~\eqref{density:sto}, with $\kappa_x\coloneq\beta\sqrt{\gamma(x)\, |\log{|\log{x}||}}$, violates condition \eqref{condition:mono} almost surely.~On the other hand, we obtain the next proposition, by observing that $\psi(\lambda,\mu)$ can be estimated in terms of~$\xi^{\lambda}$.

\begin{prop}\label{prop:twothird}
Suppose that $S$ satisfies a local functional LIL. Then, almost surely, the density $f$ in \eqref{density:sto}, with $\kappa_x \coloneq \beta\sqrt{\gamma(x)\,|\log|\log{x}||}$, adheres to condition \eqref{condition}.
\end{prop}

We start the proof of Proposition \ref{prop:twothird} with a technical lemma.

\begin{lem}\label{lem:twothird-prep}
In the context of Proposition \ref{prop:twothird}
      \begin{enumerate}[label=(\roman*)]
      % \item 
      %     $\lim_{\lambda\to 0}\sup_{\mu\in[0, 1]}\left[\int_{\f{\mu}{\mu+1}}^1\left|\f{\omega(\lambda(\mu+1) x)}{\omega(\lambda(\mu+1))} - \Gamma(x, x)\right|^3dx\right] = 0$.
      \item 
          $\lim_{\lambda\downarrow 0}\,\sup_{\mu\in[0, 1]}\,\int_{\f{\mu}{2}}^{\f{\mu+1}{2}}\Big|\f{\kappa_{2\lambda x}}
          {\kappa_{\lambda(\mu+1)}\sqrt{\gamma(x)}} - \sqrt{\frac{2^{\alpha_2}}{(\mu+1)^{\alpha_2}}}\Big|^2\,\mathrm{d}x = 0$.
      \item There exists an $\eta>0$ such that
        \begin{equation}
          \limsup_{\lambda\downarrow 0}\,\sup_{\mu\in [0, 1]}\,\int_{\f{\mu}{2}}^{\f{\mu+1}{2}} |\xi^\lambda_x| - \sqrt{\gamma(x)}\,\mathrm{d}x\le - \eta.
        \end{equation}
        \item
          $\limsup_{\lambda \downarrow 0}\,\int_0^1\big||\xi^\lambda_x|- \sqrt{\gamma(x)}\big|^{2}\,\mathrm{d}x \le 2$.
      \end{enumerate}
    \end{lem}

\smallskip

\noindent\textbf{Proof of Lemma \ref{lem:twothird-prep}.} For all small enough $\lambda>0$, all $\mu\in [0, 1]$, and all $x\in (0, 1]$, 
\begin{equation}
\begin{split}
\bigg|\f{\kappa_{2\lambda x}}{\kappa_{\lambda(\mu+1)}\sqrt{\gamma(x)}} - \sqrt{\frac{2^{\alpha_2}}{(\mu+1)^{\alpha_2}}}\bigg|
&=  \sqrt{\frac{2^{\alpha_2}}{(\mu+1)^{\alpha_2}}}\bigg| \bigg|1 + \f{\log\big|\f{\log{2\lambda x}}{\log{\lambda(\mu+1)}}\big|}
{\log|\log{\lambda(\mu+1)|}}\bigg|^{1/2}-1\bigg| \\
&\le\sqrt{2^{\alpha_2}}\bigg|\f{\log\big|\f{\log{2\lambda x}}{\log{\lambda(\mu+1)}}\big|}
{\log{|\log{\lambda(\mu+1)}|}}\bigg|^{1/2},
\end{split}
\end{equation}
where we used that $|\sqrt{|1+a|} - 1|\le \sqrt{|a|}$ for all $a\in\R$. If further $x \le \frac{\mu + 1}{2}$, then $ \frac{\log{2x/(\mu+1)}}{\log\lambda(\mu+1)} \ge 0$, and therefore
\begin{equation}
    \begin{split}
      \left|\log\left|\f{\log{2\lambda x}}{\log{\lambda(\mu+1)}}\right|\right| &= \left|\log\left|1 + \f{\log{2x/(\mu+1)}}{\log{\lambda(\mu+1)}}\right|\right|\\
      &= \log\left(1 + \left|\f{\log{2x/(\mu+1)}}{\log{\lambda(\mu+1)}}\right| \right)\\
      &\le  \left|\f{\log{2x/(\mu+1)}}{\log{\lambda(\mu+1)}}\right|.
      % &\le  \log\left(1 + \left|\f{\log{2x} - \log{(\mu+1)}}{\log{\lambda(\mu+1)}}\right|\right)\\
      % &\le \f{|\log{2x}| + |\log{(\mu+1)|}}{|\log{\lambda(\mu+1)}|}
    \end{split}
\end{equation}
Thus, we deduce
\begin{equation}
\bigg|\f{\kappa_{2\lambda x}}{\kappa_{\lambda(\mu+1)}\sqrt{\gamma(x)}} - \sqrt{\frac{2^{\alpha_2}}{(\mu+1)^{\alpha_2}}}\bigg| 
\le \sqrt{2^{\alpha_2}}\f{(|\log{2x}| + |\log{(\mu+1)}|)^{1/2}}{|\log{\lambda(\mu+1)}|^{1/2}\,|\log{|\log{\lambda(\mu+1)}|}|^{1/2}}.
\end{equation}
Result (i) follows immediately.

\medskip

The functional 
\begin{equation}\label{lem:twothird-func}
{\mathcal I}(\phi) \coloneq 
\sup_{\mu\in [0, 1]}\,\int_{\f{\mu}{2}}^{\f{\mu+1}{2}} |\phi(x)|-\sqrt{\gamma(x)}\,\mathrm{d}x
\end{equation}
on $C([0,1])$ is continuous, so that the local functional LIL implies
\begin{equation}
\limsup_{\lambda\downarrow 0}\,\sup_{\mu\in [0,1]}\,\int_{\f{\mu}{2}}^{\f{\mu+1}{2}} |\xi^\lambda_x| - \sqrt{\gamma(x)}\,\mathrm{d}x 
=\sup_{\phi\in K}\,\sup_{\mu\in [0, 1]}\,\int_{\f{\mu}{2}}^{\f{\mu+1}{2}} |\phi(x)| - \sqrt{\gamma(x)}\,\mathrm{d}x.
      \end{equation}
Notice that for all $\phi\in K$,
\begin{equation}\label{lem:twothird-cs}
|\phi(x)| = |\langle\phi, \Gamma(x, \cdot)\rangle|
\le \langle\phi,\phi\rangle^{1/2}\cdot\langle\Gamma(x, \cdot),\Gamma(x, \cdot)\rangle^{1/2}
\le\sqrt{\gamma(x)},\quad x\in [0, 1],
\end{equation}
thanks to $\langle\Gamma(x, \cdot),\Gamma(x, \cdot)\rangle=\Gamma(x,x)=\gamma(x)$.
Moreover, it is enough to prove that
\begin{equation}\label{lem:twothird-prelim-strict}
\sup_{\phi\in K}\,\sup_{a\in [0, 1/2]}\,\int_{a}^{a+1/2} |\phi(x)| - \sqrt{\gamma(x)}\,\mathrm{d}x < 0.
\end{equation}

 If the supremum in \eqref{lem:twothird-prelim-strict} was zero, then the continuity of the underlying functional on the compact $K\times[0,1/2]$ would yield the existence of some $\phi\in K$ and some  $a\in[0, 1/2]$ such that
\begin{equation}
\int_a^{a+1/2} |\phi(x)| - \sqrt{\gamma(x)}\,\mathrm{d}x = 0,
\end{equation}
and thus the Cauchy-Schwarz inequalities in \eqref{lem:twothird-cs} would hold with equality for Lebesgue almost every $x\in[a,a+1/2]$.~As a consequence, $\{\Gamma(x,\cdot)\}_x$ would be pairwise linearly dependent for these~$x$, in contradiction to the assumed non-degeneracy of $\Gamma$. This proves (ii).
 %      We see already that the right-hand side is non-positive because for all $\phi\in K$ and $x\in [0, 1]$
 %      \begin{equation}
 %        \phi(x)\le \sqrt{\int_0^xdy}\sqrt{\int_0^x\dot{\phi}(y)^2dy}\le \sqrt{x}.
 %      \end{equation}
 %      Let $\phi\in K$. To get a non-trivial upper-bound, define $\psi(x) = \int_0^x|\dot{\phi}|(y)dy$ and write for all $x_0\in [0, 1]$
 %      \begin{equation}
 %        \begin{split}
 %          \int_{x_0}^1|\phi(x)|dx\le \int_{x_0}^1\psi(y)dy &= \psi(1) - x_0\psi(x_0) - \int_{x_0}^1y|\dot{\phi}(y)|dy\\
 %                                                           &= \int_{x_0}^1(1-y)|\dot{\phi}(y)|dy + \int_0^{x_0}(1-x_0)|\dot{\phi}(y)|dy\\
 %                                                           &= \int_{0}^1\min(1-y, 1- x_0)|\dot{\phi}(y)|dy\\
 %                                                           &\le \sqrt{\int_{0}^1\min(1-y, 1- x_0)^2dy}\\
 %                                                           & \coloneqq \Psi(x_0) = \sqrt{x_0(1-x_0)^2 + \inv{3}(1-x_0)^3}.
 %        \end{split}
 %      \end{equation}

 %      Therefore,
 %      \begin{equation}
 %        \begin{split}
 % \sup_{\phi\in K}\sup_{\mu\in [0, 1]}\left[(\mu+1)\int_{\f{\mu}{\mu+1}}^{1}(|\phi(x)| - \sqrt{x})dx\right] &\le \sup_{\mu\in [0, 1]}\left[(\mu+1)\left(\Psi\left(\f{\mu}{\mu+1}\right) - \f{2}{3}\left(1 - \left(\f{\mu}{\mu+1}\right)^{3/2}\right)\right)\right]\\
 % &\coloneqq -\eta < 0.
 %        \end{split}
 %      \end{equation}
      % This proves (ii).
      
\medskip

To obtain (iii) we apply the local functional LIL to the continuous functional
\begin{equation}
C([0, 1])\to\mathbb{R},\quad\phi\mapsto
\int_0^1 \big||\phi(x)| - \sqrt{\gamma(x)}\big|^{2}\,\mathrm{d}x,
\end{equation}
and use $|\phi(x)|\le \sqrt{\gamma(x)}$, $x\in[0,1]$ for $\phi\in K$ to easily get
      \begin{equation}
\int_0^1 \big||\phi(x)| - \sqrt{\gamma(x)}\big|^{2}\,\mathrm{d}x \le 2\sup_{x\in [0, 1]} \gamma(x)
= 2
      \end{equation}
for all those $\phi$. \qed 

\medskip

We are now ready for the proof of Proposition~\ref{prop:twothird}.
    
\medskip

\noindent\textbf{Proof of Proposition~\ref{prop:twothird}.} We only need to show \eqref{condition:mean}.~Throughout the proof we take $\lambda_0>0$ to be small enough so that $|\log{\lambda x}| \ge 1$, $\lambda\in[0,\lambda_0)$, $x\in[0,2]$; $\kappa$ is non-decreasing on $[0,2\lambda_0]$; and $f(x)\le 1+ S_x-\kappa_x$, $x\in[0,2\lambda_0]$. Then,
\begin{equation}
\kappa_{\lambda x} = \beta\sqrt{\gamma(\lambda x)\,\log{|\log{\lambda x}}|} 
\ge \sqrt{\gamma(x)}\,\kappa_\lambda\left|\f{\log{|-\log{\lambda} - \log{2}|}}{\log|\log{\lambda}|}\right|^{1/2}
= \sqrt{\gamma(x)}\,\kappa_\lambda\,q_\lambda,
\end{equation}
where 
\begin{equation}\label{prop-q}
q_\lambda \coloneq 
\bigg|1 + \f{\log{\big|1 + \f{\log 2}{\log{\lambda}}\big|}}{\log{|\log{\lambda}|}}\bigg|^{1/2} \underset{\lambda\downarrow0}{\longrightarrow} 1. 
\end{equation}
It follows that, for $\lambda\in[0,\lambda_0)$ and $x\in[0,2]$,
\begin{equation}
f(\lambda x) \le 1 + \kappa_{\lambda x}\,
\left(\f{|S_{\lambda x}|}{\kappa_{\lambda x}} - 1\right)                         
\le 1 + \f{\kappa_{\lambda x}}{q_\lambda}\,\bigg(\f{|\xi^\lambda_x|}{\sqrt{\gamma(x)}} 
- q_\lambda\bigg).
\end{equation}

\smallskip

Let $\widehat{\xi}^{\lambda}_x := 2^{-\alpha_2/2}\xi^\lambda_{2x}$. Then, for $\mu\in[0,1]$,
\begin{equation}\label{eq for rmk}
\begin{split}
\psi(\lambda, \mu) - 1&\le \inv{q_\lambda}\,\int_\mu^{\mu + 1} \kappa_{\lambda x}
\,\bigg(\f{|\xi^\lambda_x|}{\sqrt{\gamma(x)}} 
- q_\lambda\bigg)\,\mathrm{d}x \\
&= \f{2}{q_\lambda}\,\int_{\f{\mu}{2}}^{\frac{\mu+1}{2}} \kappa_{2\lambda x}\, \bigg(\f{|\xi^\lambda_{2x}|}{\sqrt{\gamma(2x)}} - q_\lambda\bigg)\,\mathrm{d}x \\
&= \f{2}{q_\lambda}\,\int_{\f{\mu}{2}}^{\frac{\mu+1}{2}} \kappa_{2\lambda x}\, \bigg(\f{|\widehat{\xi}^\lambda_{x}|}{\sqrt{\gamma(x)}} - q_\lambda\bigg)\,\mathrm{d}x \\
&\le \f{2}{q_\lambda}\,\int_{\f{\mu}{2}}^{\frac{\mu+1}{2}} \kappa_{2\lambda x}\,
\bigg(\f{|\widehat{\xi}^\lambda_{x}|}{\sqrt{\gamma(x)}} - 1\bigg)\,\mathrm{d}x 
+ \left|1 - \inv{q_\lambda}\right|\,\kappa_{\lambda(\mu+1)}. 
\end{split}
\end{equation}
Next, we abbreviate $2^{\alpha_2/2}/(\mu+1)^{\alpha_2/2}$ by $\zeta(\mu)$ and rewrite 
\begin{equation}
\begin{split}
&\;\int_{\f{\mu}{2}}^{\frac{\mu+1}{2}} \kappa_{2\lambda x}\,
\bigg(\f{|\widehat{\xi}^\lambda_{x}|}{\sqrt{\gamma(x)}} - 1\bigg)\,\mathrm{d}x \\
&= \zeta(\mu)\,\kappa_{\lambda(\mu + 1)}\,\bigg(\int_{\f{\mu}{2}}^{\frac{\mu+1}{2}} |\widehat{\xi}^\lambda_x| 
- \sqrt{\gamma(x)}\,\mathrm{d}x
+ \inv{\zeta(\mu)}\int_{\f{\mu}{2}}^{\f{\mu+1}{2}} \bigg(\f{\kappa_{2\lambda x}}
{\kappa_{\lambda(\mu+1)}\sqrt{\gamma(x)}}-\zeta(\mu)\!\bigg)
\big(|\widehat{\xi}^\lambda_{x}| - \sqrt{\gamma(x)}\big)\,\mathrm{d}x\bigg).
\end{split}
\end{equation}
        % Estimate
        % \begin{equation}
        %   \begin{split}
        %     \left|\int_{\f{\mu}{\mu+1}}^{1}\left(\f{\omega(\lambda(\mu+1) x)}{\omega(\lambda(\mu+1))}-x^H\right)\left(\f{|\xi^\lambda_{x}|}{\Gamma(x, x)} - 1\right)dx\right| &\le \left(\int_{\f{\mu}{\mu+1}}^1\left|\f{\omega(\lambda (\mu+1)x)}{\omega(\lambda(\mu+1))} - \Gamma(x, x)\right|^3dx\right)^{1/3}\\
        %     &\times\left(\int_0^1\left|\f{|\xi^\lambda_x|}{\Gamma(x, x)} - 1\right|^{3/2}dx\right)^{2/3}.
        %   \end{split}
        % \end{equation}
        
\smallskip        
        
Using that $1\le \zeta(\mu) \le 2^{\alpha_2/2} $, that $(\widehat{\xi}^{\lambda}_x)_{x\in [0,1],\lambda\in (0, 1]}\stackrel{d}{=}(\xi^\lambda_x)_{x\in [0,1],\lambda\in (0, 1]}$ by the scaling relation~\eqref{scaling-prop}, and the Cauchy-Schwarz inequality in conjunction with Lemma~\ref{lem:twothird-prep}(i),(iii) we obtain
\begin{equation}\label{prop:twothird-eureka}
\limsup_{\lambda\downarrow 0} \sup_{\mu\in[0, 1]}\,
\bigg|\inv{\zeta(\mu)}\int_{\f{\mu}{2}}^{\frac{\mu+1}{2}}
\bigg(\f{\kappa_{2\lambda x}}
{\kappa_{\lambda(\mu+1)}\sqrt{\gamma(x)}}-\zeta(\mu)\!\bigg)
\big(|\widehat{\xi}^\lambda_{x}| - \sqrt{\gamma(x)}\big)\,\mathrm{d}x \bigg|= 0.
\end{equation}
In conclusion, 
\begin{equation}
\begin{split}
 \f{\psi(\lambda, \mu) - 1}{2^{\alpha_2/2}\,\kappa_{\lambda(\mu+1)}} &\le \f{\psi(\lambda, \mu) - 1}{\zeta(\mu)\,\kappa_{\lambda(\mu+1)}} \\
&\le \frac{2}{q_\lambda}\,
\sup_{\mu\in [0, 1]}\,\int_{\f{\mu}{2}}^{\frac{\mu+1}{2}}
|\widehat{\xi}_x^\lambda| - \sqrt{\gamma(x)}\,\mathrm{d}x \\
&\quad\, + \f{2}{q_\lambda}\,\sup_{\mu\in[0, 1]}\,
\bigg|\inv{\zeta(\mu)}\int_{\f{\mu}{2}}^{\frac{\mu+1}{2}}
\bigg(\f{\kappa_{2\lambda x}}
{\kappa_{\lambda(\mu+1)}\sqrt{\gamma(x)}}\!-\!\zeta(\mu)\!\bigg)
\big(|\widehat{\xi}^\lambda_{x}| \!-\! \sqrt{\gamma(x)}\big)\,\mathrm{d}x \bigg| 
+  \left|1 - \inv{q_\lambda}\right|,
\end{split}
\end{equation}
for which \eqref{prop-q}, $(\widehat{\xi}^{\lambda}_x)_{x\in [0,1],\lambda\in (0, 1]}\stackrel{d}{=}(\xi^\lambda_x)_{x\in [0,1],\lambda\in (0, 1]}$, Lemma \ref{lem:twothird-prep}(ii) and \eqref{prop:twothird-eureka} yield the existence of a $\lambda_0>0$ such that
\begin{equation}
\psi(\lambda, \mu) \le 1 - 2^{\alpha_2/2}\,\kappa_{\lambda(\mu+1)}\,\frac{\eta}{2}
\end{equation}
for all $\lambda\in [0, \lambda_0)$ and $\mu\in [0, 1]$. \qed 

\begin{rmk}\label{rmk:other-proc}
\begin{enumerate}[wide,label=(\alph*),labelindent=0pt]
\item Our proof of Proposition \ref{prop:twothird} also applies to densities $f$ with the property 
\begin{equation}
 f(x) = \bigg(1 + \widetilde{\kappa}_x\,\bigg(\f{S_x}{\kappa_x} - 1\bigg)\!\bigg)_+\wedge 1,\quad x\in[0,1],
\end{equation}
for a non-negative non-decreasing function $\widetilde{\kappa}$ obeying
\begin{equation}
\lim_{\lambda\downarrow 0}\,\sup_{\mu\in[0, 1]}\,
\int_{\f{\mu}{2}}^{\f{\mu+1}{2}}\bigg|\f{\widetilde{\kappa}_{2\lambda x}}
{\widetilde{\kappa}_{\lambda(\mu+1)}\sqrt{\gamma(x)}} - \sqrt{\f{2^{\alpha_2}}{(\mu+1)^{\alpha_2}}}\bigg|^2\,\mathrm{d}x = 0.
\end{equation}
In addition, one can cover densities $f$ with
\begin{equation}
f(x) = \f{|S_x|}{\kappa_x}\wedge 1,\quad x\in[0, 1]
\end{equation}
by repeating the proof of Lemma \ref{lem:twothird-prep}(ii) for the final line in \eqref{eq for rmk} with $1$ in place of $\kappa$. \\
\item By using a very similar method, we can verify condition \eqref{condition} for densities $f$ such that
\begin{equation}
f(x) = \f{|S_{1/x}|}{\kappa_{1/x}}\wedge 1,\quad x\in (0, 1],
\end{equation}
where $S$ satisfies the scaling property \eqref{scaling-prop} and a local functional LIL ``at infinity'': Almost surely, the family $\big\{\big(\f{S_{r x}}{\kappa_r}\big)_{x\in[0,1]}\big\}_{r\ge 3}$ is relatively compact in $C([0, 1])$ with the set of limit points $K$ as above.~The local functional LIL at infinity is known for various classes of Gaussian processes $S$, including fractional Brownian motion (see~\cite[Example 4.36]{malyarenko2012invariant}), semi-stable Gaussian processes (see~\cite[Theorem 4]{oodaira1972strassen}), Gaussian processes that are not necessarily semi-stable (see~\cite[Theorem 4]{oodaira1973law}) but for which \cite[Condition (A-1)]{oodaira1973law} makes an adaptation of our proof possible, and rescalings of Brownian motion (see \cite[Theorems 1--3]{bulinski1980new}). \\
\item Another interesting process admitting a local functional LIL at infinity is iterated Brownian motion (see~\cite[Theorem 1.1]{hu1995laws}).~In this case, our proof can be adjusted as follows.~Let $(W^1_x)_{x\in \R}$ and $(W^2_x)_{x\ge 0}$ be two independent standard Brownian motions.~Define
\begin{equation}
S_x = W^1_{W^2_x},\quad x\ge 0 \qquad\qquad\qquad\;\;\;
\end{equation}
and 
\begin{equation}
\kappa_x = 2^{3/4} x^{1/4}(\log\log x)^{3/4},\quad x\ge3.
\end{equation}
The relevant compact subset $K$ of $C([0, 1])$ is then given by 
\begin{equation}
K = \bigg\{f \circ g:\; f\in C([-1, 1]),\, g\in C([0, 1]),\, f(0) = 0,\, g(0) = 0,\,
\int_{-1}^{1}f'(x)^2\,\mathrm{d}x + \int_0^1 g'(x)^2\,\mathrm{d}x \le 1\bigg\}.
\end{equation}
Indeed, \cite[Theorem 1.1]{hu1995laws} implies that, almost surely,
\begin{equation}
\limsup_{r\to\infty}\,\mathcal{I}\left(\f{S_{r\cdot}}{\kappa_r}\right) 
= \sup_{\phi\in K}\,\mathcal{I}(\phi),
\end{equation}
for any continuous functional $\mathcal{I}\!:C([0, 1]) \to \R$.~This allows us to redo the proofs of Lemma~\ref{lem:twothird-prep} and Proposition~\ref{prop:twothird}. In particular, the inequalities in \eqref{lem:twothird-cs} can be replaced by
\begin{equation}
\begin{split}
\phi(x)\!=\!\int_0^1 f'(y)\,\mathbf{1}_{\{y\le g(x)\}}\,\mathrm{d}y
\le\! \bigg(\int_0^1 f'(y)^2\,\mathrm{d}y\bigg)^{1/2} \sqrt{g(x)}
\le\! \sqrt{\bigg(\int_0^x g'(y)^2\,\mathrm{d}y\bigg)^{1/2} x^{1/2}}
\le\! x^{1/4}, 
\; x\in[0,1],
\end{split}
\end{equation}
for all $\phi = f\circ g\in K$.
\end{enumerate}        
\end{rmk}

\subsection{Initial conditions constructed from periodic functions }\label{subsection:check-sinus}
%%%%%%%%%%%%%%%%%%%%%%%%%%%%%%%%%%%

Let $\Psi\!:[0,\infty)\to[-1,1]$ be a periodic function with $\sup_{x\ge0} \int_0^x \Psi(y)\,\mathrm{d}y<\infty$ and $\limsup_{x\to\infty}\,\Psi(x) = 1$.~In this subsection, we show that, for any $\alpha>0$, the oscillating probability density given by
\begin{equation}\label{def:sin}
f(x) = \inv{2}\left(1 + \Psi\left(\inv{x^{\alpha}}\right)\!\right),\quad x\in(0,a]
\end{equation}
satisfies condition \eqref{condition}. The parameter $\alpha$ controls how fast the density oscillates (cf.~Remark \ref{rmk:oscil}).

\begin{prop}
Every probability density $f$ defined by \eqref{def:sin} obeys condition \eqref{condition}.
  %         \begin{equation}
  % \limsup_{\lambda\to 0}\left(\sup_{\mu\in [0, C]}\psi(\lambda, \mu)\right) < 1 -  \f{1-\alpha_1\alpha_2}{2(C+1)}.
  %         \end{equation}
\end{prop}
        
\noindent\textbf{Proof.} We only need to check \eqref{condition:mean}. To this end, for $\lambda\in\big(0,\frac{a}{2}\big)$ and $\mu\in [0, 1]$, we compute
\begin{align*}
 \psi(\lambda,\mu) - \inv{2} 
 = \frac{1}{2} \int_{\mu}^{\mu + 1} \Psi\left(\inv{\lambda^\alpha x^{\alpha}}\right)\,\mathrm{d}x 
 = \inv{2\alpha\lambda} 
 \int_{\inv{\lambda^\alpha(\mu + 1)^{\alpha}}}
 ^{\inv{\lambda^\alpha\mu^{\alpha}}}\f{\Psi(x)}{x^{\inv{\alpha} + 1}}\,\mathrm{d}x.
\end{align*}
Integrating by parts, writing $H(x)$ for $\int_0^x\Psi(y)\,\mathrm{d}y$, and using $\mu+1\le 2$ we get
\begin{align*}
\psi(\lambda,\mu) - \inv{2} &=\inv{2\alpha\lambda}\left[\lambda^{\alpha+1}\mu^{\alpha+1}H\left(\inv{\lambda^\alpha\mu^{\alpha}}\right) - \lambda^{\alpha+1}(\mu+1)^{\alpha+1}H\left(\inv{\lambda^\alpha(\mu + 1)^{\alpha}}\right) \right]\\
&\quad + \inv{2\alpha\lambda}\left(\inv{\alpha} + 1\right) 
\int_{\inv{\lambda^\alpha(\mu + 1)^\alpha}}
^{\inv{\lambda^\alpha\mu^\alpha}}\f{H(x)}{x^{\inv{\alpha}+2}}\,\mathrm{d}x \\
&\le\sup_{x\ge0} H(x)\,\frac{\lambda^\alpha}{\alpha}\,2^{\alpha+1}.
\end{align*}
Therefore, it holds 
\begin{equation}
\sup_{\mu\in [0, 1]}\psi(\lambda, \mu)< \f{3}{4}
\end{equation}
for all $\lambda\ge0$ small enough. \qed      
\section{Refined analysis for some piecewise constant initial conditions} \label{section:toy}
%%%%%%%%%%%%%%%%%%%%%%%%%%%%%%%%%%%%%%%

This section is devoted to the well-posedness question for oscillatory and piecewise constant probability densities defined by
\begin{equation}\label{density:toy'}
	f(x) = 
	\begin{cases}
		\alpha_1,\quad x\in \underset{n\ge1}{\bigcup}\,[a_{2n}, a_{2n-1}), \\
		\alpha_2,\quad x\in \underset{n\ge1}{\bigcup}\,[a_{2n+1}, a_{2n}),
	\end{cases}
\end{equation}
where $0 < \alpha_1<1<\alpha_2$, $a_{2n-1} = r^{n-1}a_1$, $a_{2n} = pr^{n-1}a_1$, and $r=pq$, $p, q\in (0, 1)$.~Such densities are of interest because they can violate both  \eqref{condition:less-one} and \eqref{condition:mean}, thus necessitating additional arguments to prove the uniqueness of the associated physical solution.~Note that the CDF $F$ is piecewise linear and oscillates between the half-lines $y = \beta_1 x$ and $y = \beta_2 x$, with $0<\beta_1< \beta_2$ given by
\begin{eqnarray}
&& \beta_1 =  \inv{1-pq}\big(\alpha_2p(1-q) + \alpha_1(1-p)\big), \\
&& \beta_2 =  \inv{1-pq}\big(\alpha_2(1-q) + \alpha_1q(1-p)\big).
\end{eqnarray}
For technical reasons (see Proposition~\ref{lem:toy-square-root-bounds} below), we assume in the following that $\beta_2<1$, namely
\begin{equation}\label{conditon-beta2}
\alpha_2 < 1 + q\f{1-p}{1-q}(1-\alpha_1).
\end{equation}
Condition \eqref{condition:less-one} is not satisfied by $f$. For $q\in(0,1/2]$, condition \eqref{condition:mean} fails for it as well.

\begin{prop}
For $q\in(0,1/2]$, the density $f$ defined by \eqref{density:toy'} violates condition \eqref{condition:mean}.
\end{prop}

\noindent\textbf{Proof.} Take $\lambda=\frac{1-q}{q}\,a_{2n+1}$ for an integer $n\ge1$ and set $\widetilde{\mu}=\f{q}{1-q}\in(0,1]$. Observe that $\lambda\widetilde{\mu}=a_{2n+1}$, whereas $\lambda(\widetilde{\mu}+1) = a_{2n+1}\big(1+\frac{1-q}{q}\big)=a_{2n}$. Thus,
\begin{equation}
\int_{\widetilde{\mu}}^{\widetilde{\mu} + 1} f(\lambda x)\,\mathrm{d}x =\alpha_2>1.
\end{equation}
Consequently, also
\begin{equation}
\sup_{\mu\in [0, 1]}\,\int_{\mu}^{\mu+1}f(\lambda x)\,\mathrm{d}x = \alpha_2>1.
\end{equation}
Hence, condition \eqref{condition:mean} cannot hold. \qed

\medskip

Nevertheless, we are able to prove Proposition~\ref{prop:toy}.~Our proof relies on the next proposition, akin to Proposition~\ref{prop:iteration}.
\begin{prop}\label{prop:iteration-toy}
For any $\alpha_2>1$ close enough to $1$,
\begin{equation}\label{def:delta0}
\sup_{t\in (0, T]}\,\sup_{h>0}\,\E{\f{F(Y_{t} + h) - F(Y_{t})}{h}}=:\delta_0 < 1.
\end{equation}
\end{prop}
\noindent Once this result is proved, the desired uniqueness on $[0,T]$ can be shown by proceeding as in the proof of Theorem \ref{thm:uniqueness}, only with $1-\delta_0$ in place of $\Phi(t,\lambda)$. The strategy of the proof of Proposition~\ref{prop:iteration-toy}, in turn, lies in finding a set $G\subset[0,\infty)$ such that for $\alpha_2>1$ close enough to $1$,
\begin{equation}\label{toy-step1}
\sup_{y\in G}\,\sup_{h>0:\,y + h \le a_1}\,\f{F(y + h) - F(y)}{h} \eqqcolon L < 1.
\end{equation}
Then, estimating the expectation in \eqref{def:delta0} according to
\begin{equation}
\E{\f{F(Y_{t} + h) - F(Y_{t})}{h}}\le \alpha_{2} - (\alpha_{2}-L)\,\Prob{Y_{t}\in G}
\end{equation}
it remains to check that $Y_t$ falls into $G$ with a sufficiently high probability, namely
\begin{equation}\label{toy-step2}
\inf_{t\in (0, T]}\,\Prob{Y_{t}\in G}> \f{\alpha_{2}-1}{\alpha_{2}-L}.
\end{equation}
The two assertions \eqref{toy-step1} and \eqref{toy-step2} are the subjects of Subsections \ref{sec:4.1} and \ref{sec:4.2}, respectively. 

%%%%%%%%%%%%%%%%%%%%%
\subsection{Proof of \eqref{toy-step1}}\label{sec:4.1}
%%%%%%%%%%%%%%%%%%%%%

\begin{lem}\label{lem:rapport}
Let $G = \bigcup_{n\ge 1} [a_{2n+2},\varrho a_{2n+1}]\cup[a_2,\infty)$, where $\varrho\coloneq \f{1 + p}{2}$. Then,
\begin{equation}\label{supG}
\sup_{y\in G}\,\sup\limits_{h>0:\,y+h\le a_1}\,
\f{F(y + h) - F(y)}{h} = \f{(1-q)\alpha_2 + q(1-\varrho)\alpha_1}{1-q\varrho} \eqqcolon L.
\end{equation}
Moreover, for $\alpha_2>1$ close enough to $1$, it holds $L<1$.
\end{lem}

\noindent\textbf{Proof.}~It suffices to show \eqref{supG} with $G\backslash[a_2,\infty)$ in place of $G$.~To this end, fix an $n\ge 1$ and a $y\in [a_{2n+2}, \varrho a_{2n+1}]$. Define the function 
\begin{equation}
\theta:\,(0,a_1-y]\to[0,\infty),\quad h\mapsto\f{F(y+h) - F(y)}{h}.
\end{equation}
By the definition of $F$, we have for $k=1,\,2,\,\ldots,\,n$:
\begin{equation}
\begin{cases}
\theta'(h) \ge 0,\quad y+h\in (a_{2k+1}, a_{2k}), \\
\theta'(h) \le 0,\quad y+h\in (a_{2k}, a_{2k-1}).
\end{cases}
\end{equation}
Therefore, 
\begin{equation}
\sup\limits_{h>0:\,y+h\le a_1}\,\f{F(y + h) - F(y)}{h} 
= \sup_{a_{2k}\ge a_{2n}}\f{F(a_{2k}) - F(y)}{a_{2k}-y}.
\end{equation}

\smallskip

Notice now that the sequence
  \begin{equation}
    \left(\f{F(a_{2k}) - F(y)}{a_{2k}-y}\right)_{k=1,\,2,\,\ldots,\,n}
  \end{equation}
  is non-decreasing. Indeed, for $k=2,\,3,\,\ldots,\,n$,
  \begin{equation}
        \f{F(a_{2k}) - F(y)}{a_{2k}-y}  - \f{F(a_{2k-2}) - F(y)}{a_{2k-2}-y} = \f{(y\beta_2-F(y))(a_{2k-2} - a_{2k})}{(a_{2k} - y)(a_{2k-2}-y)}\ge 0.
  \end{equation}
We conclude
  \begin{equation}
    \sup\limits_{h>0:\,y+h\le a_1}\,\f{F(y + h) - F(y)}{h} = \f{F(a_{2n}) - F(y)}{a_{2n}-y}.
  \end{equation}
Since the right-hand side is non-decreasing in $y$ on $[a_{2n+2}, \varrho a_{2n+1}]$,
  \begin{equation}
    \sup_{y\in [a_{2n+2}, \varrho a_{2n+1}]}\,\sup\limits_{h>0:\,y+h\le a_1}\,\f{F(y + h) - F(y)}{h}  
    = \f{F(a_{2n}) - F(\varrho a_{2n+1})}{a_{2n} - \varrho a_{2n+1}} 
    = \f{(1-q)\alpha_2 +q(1-\varrho)\alpha_1}{1-q\varrho}.
  \end{equation}
This proves the first statement. The second one is straightforward to verify. \qed

%%%%%%%%%%%%%%%%%%%%
\subsection{Proof of (\ref{toy-step2})}\label{sec:4.2}
%%%%%%%%%%%%%%%%%%%%

The key step in deriving \eqref{toy-step2} is an estimate of the probabilities
\begin{equation}\label{density-estimate}
\PP\big(Y_t\in [a\sqrt{t}, b\sqrt{t}]\big),\quad 0<a<b,\quad t\in(0,T].
\end{equation}
For that purpose, we establish the $1/2$--H\"older continuity of the frontier $\Lambda_{}$ on $[0,T]$. 
%   Let $T > 0$. We prove in this section the existence of two constants $c_1$ and $c_2$ such that for all $t\in [0, T]$. 
%   \begin{equation}
%     c_1\sqrt{t}\le \Lambda_{t}\le c_2\sqrt{t}.
%   \end{equation}
% Initial density (\ref{density:toy}) leads to a CDF oscillating between $y = \beta_1x$ and $y = \beta_2x$ for $\beta_1$ and $\beta_2$ given by
% \begin{equation}
%   \beta_1 = \inv{1-pq}\left[\alpha_2p(1-q) + \alpha_1(1-p)\right],
% \end{equation}
% and
% \begin{equation}
%   \beta_2 = \inv{1-pq}\left[\alpha_2(1-q) + \alpha_1q(1-p)\right].
% \end{equation}
% The existence of $c_1$ is therefore immediate by Lemma~\ref{lem:square-root-lb}. The derivation of the square-root upper-bound requires the technical assumption $\beta_1<1$, equivalent to
% \begin{equation}
%   \label{cond:beta2}
%   \alpha_2 < 1 + q\f{1-p}{1-q}(1-\alpha_1).
% \end{equation}
As a preparation for the latter, we introduce for each $t\in[0,T]$ the function
\begin{equation}
F_t:\;[0,\infty)\to[0,1],\quad x\mapsto
\Prob{0<X_t\le x} = \E{F(\Lambda_{t} - B_t + x) - F(\Lambda_{t} -B_t)}
\end{equation}
and notice immediately that $F_t'(x)\le\alpha_2$. Moreover, we have the following bound. 

\begin{lem}\label{lem:square-root-ub}
For all $t\in[0,T]$, it holds 
  \begin{equation}
    \Lambda_{t+h} - \Lambda_{t} - F_t(\Lambda_{t+h} - \Lambda_{t})\le \alpha_2\sqrt{\f{2}{\pi}}\,\sqrt{h},\quad h>0.
  \end{equation}
\end{lem}

\noindent\textbf{Proof.} We start with the inequalities
\begin{equation}
\begin{split}
\Lambda_{t+h} - \Lambda_{t} 
&= \PP\big(\sup_{0\le s\le t}(-B_s + \Lambda_{s}) < \Xzero\le \sup_{0\le s\le t+ h} (-B_s + \Lambda_{s})\big) \\
&= \PP\big(\sup_{0\le s\le t}(-B_s + \Lambda_{s})+B_t-\Lambda_{t} < \Xzero+B_t-\Lambda_{t}
\le \sup_{t\le s\le t+ h}(-B_s + \Lambda_{s})+B_t-\Lambda_{t}\big) \\
&\le \PP\big(0 < X_t\le \sup_{t\le s\le t+ h}\left(B_t-B_s\right) + \Lambda_{t+h} - \Lambda_{t}\big) \\
&= F_t(\Lambda_{t+h} - \Lambda_{t}) + \PP\big(\{\Lambda_{t+h} - \Lambda_{t}<X_t\} \cap \big\{X_t - (\Lambda_{t+h}-\Lambda_{t}) \le  
\sup_{t\le s\le t+ h} (B_t-B_s)\big\}\big).
\end{split}
\end{equation}
Consequently,
\begin{equation}
\Lambda_{t+h} - \Lambda_{t} - F_t(\Lambda_{t+h} - \Lambda_{t})
\le \int_0^\infty \PP\big(x\le \sup_{t\le s\le t+h}(B_t-B_s)\big)\,\mathrm{d}F_t(x + \Lambda_{t+h} - \Lambda_{t})
\le \alpha_2\sqrt{\f{2}{\pi}}\,\sqrt{h},
\end{equation}
as stated in the lemma. \qed

\medskip

As a direct implication, we obtain the square root behavior of the frontier $\Lambda_{}$.

\begin{prop}\label{lem:toy-square-root-bounds} 
For any $\alpha_2>1$ close enough to $1$, there exist $0<c_1\le c_2<\infty$ such that 
\begin{equation}
c_1\sqrt{t}\le \Lambda_{t}\le c_2\sqrt{t},\quad t\in[0,T].
\end{equation}
\end{prop}

\noindent\textbf{Proof.} For the lower bound, we notice that $Z_t\ge \sup_{0\le s\le t}(-B_s)$, and hence,
\begin{equation}
\Lambda_{t} = \E{F(Z_t)}\ge \EE\big[F\big(\sup_{0\le s\le t}(-B_s)\big)\big]
\ge\beta_1\EE\big[\sup_{0\le s\le t}(-B_s)\big]
=\beta_1\sqrt{\f{2}{\pi}}\,\sqrt{t},\quad t\in[0,T].
\end{equation}
For the upper bound, we apply Lemma~\ref{lem:square-root-ub} with $t=0$ and get 
\begin{equation}\label{def:c2}
\Lambda_h \le \f{\alpha_2}{1-\beta_2}\sqrttp\,\sqrt{h},\quad h\in(0,T]
\end{equation}
thanks to $\Lambda_{0}=0$ and $F_0(x)=F(x)\le\beta_2 x$. \qed

\medskip  

The $1/2$--H\"older continuity of $\Lambda_{}$ on $[0,T]$ is deduced in a similar way from the next proposition. 

\begin{prop}\label{lem:slope}
For any $\alpha_2>1$ close enough to $1$, there exists a $\beta\in[0, 1)$ such that
\begin{enumerate}[label=(\roman*)]
\item $F_t(x) \le \beta x$, $x\ge 0$, $t\in[0,T]$, and
\item $\E{f(\Lminus{t} + B_t)} \le \beta$, $t\in(0, T]$.
\end{enumerate}
In particular, $\Lambda_{}$ is continuous on $[0,T]$.
\end{prop}

\noindent\textbf{Proof.} Fix a $C\in(0,\infty)$ and consider a $t\in(0,T]$. Then, for $x>C\Lambda_{t}$,
\begin{equation}
F_t(x)\le \E{F(\Lambda_{t} - B_t + x)} 
\le \beta_2\E{(\Lambda_{t} - B_t + x)_+} 
\le \beta_2\bigg(\f{1+C}{C}x + \f{1}{\sqrt{2\pi}}\,\sqrt{t}\bigg).
\end{equation}
In view of the square root lower bound $\sqrt{t}\le\f{\Lambda_{t}}{c_1}$, we have for $x>C\Lambda_{t}$,
\begin{equation}
F_t(x)\le \beta_2\bigg(\f{1+C}{C} + \inv{Cc_1\sqrt{2\pi}}\bigg)x.
\end{equation}
Since $c_1 = \beta_1\sqrt{2/\pi}\ge \alpha_1\sqrt{2/\pi}$, we conclude
\begin{equation}\label{lem:slope-case-1}
F_t(x)\le \beta_2\bigg(\f{1+C}{C} + \inv{2\alpha_1C}\bigg)x,\quad x>C\Lambda_{t}.
\end{equation}

\smallskip

Next, take $x\le C\Lambda_{t}$. By definition,
\begin{equation}\label{small x}
F_t(x) = \E{\int_{\Lambda_{t}-B_t}^{\Lambda_{t}-B_t+x} f(y)\,\mathrm{d}y}
=\E{\int_{\Lambda_{t}}^{\Lambda_{t}+x} f(y-B_t)\,\mathrm{d}y}
=\int_{\Lambda_{t}}^{\Lambda_{t}+x}\E{f(y + B_t)}\,\mathrm{d}y.
\end{equation}
Thus, it suffices to show that for any $\alpha_2>1$ close enough to $1$, there exists a $\beta\in [0, 1)$ such that 
\begin{equation}
\E{f(y+B_t)} \le \beta,\quad y\in [\Lambda_{t}, (1+C)\Lambda_{t}].
\end{equation}
Set $H = \bigcup\limits_{k\ge 1}[a_{2k}, a_{2k-1})\cup[a_1,\infty)$ and estimate $\E{f(y+B_t)}$ according to 
\begin{equation}\label{H decomp}
\E{f(y + B_t)} =\E{(f\,\mathbf{1}_H)(y+B_t)} + \E{(f\,\mathbf{1}_{H^c})(y+B_t)} 
\le \alpha_2 - (\alpha_2 - \alpha_1)\,\Prob{y + B_t\in H}\!.
\end{equation}
Our goal now is to lower bound $\Prob{y + B_t\in H}$ for $y\in [\Lambda_{t}, (1+C)\Lambda_{t}]$. We distinguish four cases.

\medskip

\noindent\textbf{Case 1:} $y\in \big[a_{2n+2}, \f{a_{2n+2}+a_{2n+1}}{2}\big)$ for some $n\ge 0$. In this case, we find
\begin{equation}\label{prop-slope-first-eq}
\Prob{y + B_t\in H}\ge 
\PP\big(B_t\in [a_{2n+2} - y, a_{2n+1} - y)\big)
\ge \PP\Big(B_t\in \Big[0, \f{a_{2n+1} - a_{2n+2}}{2}\Big)\!\Big).
\end{equation}
In view of
\begin{equation}
\f{a_{2n+1} - a_{2n+2}}{2} = \f{1-p}{2}a_{2n+1}\ge \f{1-p}{2}y\ge \f{1-p}{2}\Lambda_{t},
\end{equation}
we get
\begin{equation}\label{prop:slope-case-1}
\Prob{y + B_t\in H}\ge\PP\Big(B_t\in \Big[0, \f{1-p}{2}\Lambda_{t}\Big)\!\Big).
\end{equation}
  
\smallskip
  
\noindent\textbf{Case 2:} $y\in\big[\f{a_{2n+2}+a_{2n+1}}{2}, a_{2n+1}\big)$ for some $n\ge 0$. Similarly to the previous case, we have
\begin{equation}\label{prop:slope-case-2}
\Prob{y + B_t\in H}\ge \PP\Big(B_t\in \Big[-\f{a_{2n+1}-a_{2n+2}}{2} ,0\Big]\Big)
\ge \PP\Big(B_t\in \Big[0, \f{1-p}{2}\Lambda_{t}\Big]\Big).
\end{equation}
  
\smallskip

\noindent\textbf{Case 3:} $y\in \left[a_{2n+1}, a_{2n}\right)$ for some $n\ge 1$. In this case, 
\begin{equation}
\begin{split}
\Prob{y + B_t\in H}&\ge \Prob{B_t\in [a_{2n} - y, a_{2n-1} - y]} \\
&\ge \f{a_{2n-1} - a_{2n}}{\sqrt{2\pi t}}\,e^{-\f{(a_{2n-1} - y)^2}{2t}} \\
&\ge \f{a_{2n-1} - a_{2n}}{\sqrt{2\pi t}}\,e^{-\f{(a_{2n-1} - a_{2n+1})^2}{2t}}.
\end{split}
\end{equation}
Using 
\begin{equation}
a_{2n-1} - a_{2n}= \f{1-p}{p}a_{2n}\ge \f{1-p}{p}y\ge \f{1-p}{p}\Lambda_{t}
\end{equation}
and
\begin{equation}
a_{2n-1}-a_{2n+1} = \f{1-pq}{pq}a_{2n+1}\le \f{1-pq}{pq}y\le \f{1-pq}{pq}(1+C)\Lambda_{t}
\end{equation}
we end up with 
\begin{equation}\label{prop-slope-last-eq}
\Prob{y + B_t\in H}\ge \f{1-p}{p}\f{\Lambda_{t}}{\sqrt{2\pi t}}e^{-\f{(1-pq)^2(1+C)^2}{(pq)^2}\f{\Lambda_{t}^2}{2t}}.
\end{equation}

\smallskip

\noindent\textbf{Case 4:} $y\in[a_1,\infty)$. Here,
\begin{equation}\label{case4}
\Prob{y + B_t\in H}\ge\PP\big(B_t\in [a_1-y,\infty)\big)\ge\frac{1}{2}.
\end{equation}

\smallskip

Combining \eqref{prop:slope-case-1}, \eqref{prop:slope-case-2}, \eqref{prop-slope-last-eq} and \eqref{case4}, and employing $c_1\le \f{\Lambda_{t}}{\sqrt{t}}\le c_2$, we arrive at
\begin{equation}\label{prob-G}
\Prob{y + B_t\in H}\ge \min\bigg(\PP\Big(\mathcal{N}\in \Big[0, \f{1-p}{2}c_1\Big]\Big),  \f{1-p}{p}\f{c_1}{\sqrt{2\pi}}e^{-\f{(1-pq)^2(1+C)^2}{(pq)^2}\frac{c_2^2}{2}},\frac{1}{2}\bigg).
\end{equation}
At this point, we choose $C = \frac{(2\alpha_1+1)\beta_2}{\alpha_1(1-\beta_2)}$, so that 
\begin{equation}\label{C choice}
\beta_2\bigg(\f{1+C}{C} + \inv{2\alpha_1C}\bigg) = \f{1+ \beta_2}{2} < 1.
\end{equation}
Then, the right-hand side in \eqref{prob-G} depends on $\alpha_2$ via $c_1$, $c_2$, and $C$. For $\alpha_2\!\downarrow\!1$, the values of~$\beta_1$, $\beta_2$ tend to (distinct) limits in $(0,1)$, hence $c_1$ stays bounded away from zero, and $c_2$, $C$ stay bounded away from infinity. Therefore,
\begin{equation}
\liminf_{\alpha_2\downarrow 1}\,\Prob{y + B_t\in H}\ge
\liminf_{\alpha_2\downarrow 1}\,\min\bigg(\PP\Big(\mathcal{N}\in \Big[0, \f{1-p}{2}c_1\Big]\Big),  \f{1-p}{p}\f{c_1}{\sqrt{2\pi}}e^{-\f{(1-pq)^2(1+C)^2}{(pq)^2}\frac{c_2^2}{2}},\frac{1}{2}\bigg) > 0.
  \end{equation}
Consequently, for any $\alpha_2>1$ close enough to $1$,
\begin{equation}
\Prob{y + B_t\in H} > \f{\alpha_2 - 1}{\alpha_2-\alpha_1},
\end{equation}
yielding by \eqref{H decomp} a $\beta\in[0,1)$ such that 
\begin{equation}\label{f bound beta}
\E{f(y+B_t)} \le \beta,\quad y\in [\Lambda_{t},(1+C)\Lambda_{t}].
\end{equation}
Together with \eqref{small x}, \eqref{lem:slope-case-1} and \eqref{C choice} this finishes the proof of (i). 

\medskip

Result (ii) can be obtained by noticing that
  \begin{equation}
    c_1 = c_1\sup_{0<s<t}\f{\sqrt{s}}{\sqrt{t}}\le \sup_{0<s<t}\f{\Lambda_{s}}{\sqrt{t}}= \f{\Lminus{t}}{\sqrt{t}},
  \end{equation}
and by subsequently repeating \eqref{H decomp}--\eqref{f bound beta} mutatis mutandis.~Lastly, the final statement in the proposition is immediate from (ii) and the physical jump condition \eqref{def:physical}. \qed

\medskip

Combining Lemma~\ref{lem:square-root-ub} and Proposition~\ref{lem:slope} we deduce the next proposition. 

\begin{prop}
For any $\alpha_2>1$ close enough to $1$, there exists a $c_3\in(0,\infty)$ such that
\begin{equation}
\Lambda_{t+h}-\Lambda_{t}\le c_3\sqrt{h},\quad h\in[0,T-t],\quad t\in[0,T].
\end{equation}
Moreover, $c_3$ can be chosen according to
\begin{equation}\label{def:c3}
c_{3} = \f{\alpha_2}{1 - \beta}\sqrttp.
\end{equation}
\end{prop}

We are now ready to estimate the probabilities in \eqref{density-estimate}.

\begin{lem}\label{lem:density}
Let $U \coloneq \sup_{0\le s\le 1}(B_s + c_3\sqrt{s})$. Then, for any $\alpha_2>1$ close enough to $1$,
\begin{equation}
\PP\big(\sup_{0\le s\le t}(-B_s + \Lambda_{s})\in [a\sqrt{t}, b\sqrt{t}]\big)
\ge \Prob{\left|\mathcal{N}\right|\ge a}\,\Prob{U\le b-a},\quad 0<a<b,\quad t\in(0,T].
\end{equation}
\end{lem}

\noindent\textbf{Proof.} We fix $0<a<b$, $t\in(0,T]$, and set
\begin{equation}
\tau_a = \inf\{s>0:\,B_s + \Lambda_{s} \ge a\sqrt{t}\}.
\end{equation}
Consider the representation
\begin{equation}
\PP\big(\sup_{0\le s\le t}(-B_s + \Lambda_{s})\in [a\sqrt{t}, b\sqrt{t}]\big) 
= \PP\big(\tau_a\le t,\,\sup_{\tau_a\le s\le t}(B_s + \Lambda_{s})\le b\sqrt{t}\big).
\end{equation}
By the continuity of $\Lambda$,
\begin{equation}
B_{\tau_a} + \Lambda_{\tau_a} = a\sqrt{t},
\end{equation}
and therefore writing $W$ for the Brownian motion $B_{\tau_a + \cdot} - B_{\tau_a}$ we find 
\begin{equation}
\PP\big(\sup_{0\le s\le t}(-B_s + \Lambda_{s})\in [a\sqrt{t}, b\sqrt{t}]\big)  
= \PP\big(\tau_a\le t, \sup_{0\le s\le t - \tau_a}(W_s + \Lambda_{\tau_a+s} - \Lambda_{\tau_a})\le (b -a)\sqrt{t}\big).
\end{equation}

\smallskip

Next, we use $\Lambda_{\tau_a+s} - \Lambda_{\tau_a}\le c_3\sqrt{s}$ to deduce
\begin{equation}
\PP\big(\sup_{0\le s\le t}(-B_s + \Lambda_{s})\in [a\sqrt{t}, b\sqrt{t}]\big)  
\ge \PP\big(\tau_a\le t, \sup_{0\le s\le t}(W_s + c_3\sqrt{s})\le(b - a) \sqrt{t}\big).
\end{equation}
The trivial lower bound $\Lambda_{}\ge 0$ implies
\begin{equation}
\Prob{\tau_a\le t}= \PP\big(\sup_{0\le s\le t}\left(B_s+\Lambda_{s}\right)\ge a\sqrt{t}\big)
\ge \PP\big(\sup_{0\le s\le t}B_s\ge a\sqrt{t}\big) 
= \Prob{\left|\mathcal{N}\right|\ge a}\!.
\end{equation}
This and the independence of $W$ from $\tau_{a}$ yield
\begin{equation}
\PP\big(\sup_{0\le s\le t}(-B_s + \Lambda_{s})\in [a\sqrt{t}, b\sqrt{t}]\big)   
\ge \Prob{\tau_a\le t}\,\Prob{U\le b-a}
\ge \Prob{\left|\mathcal{N}\right|\ge a}\,\Prob{U\le b - a}\!,
\end{equation}
finishing the proof. \qed

\medskip

We conclude with the proof of \eqref{toy-step2}.

\begin{lem}\label{lem:prob-G}
For any $\alpha_2>1$ close enough to $1$, there exists a $\delta<1$ such that
\begin{equation}\label{eq: final lemma}
\inf_{t\in(0, T]}\,\Prob{Y_{t}\in G}\ge \f{\alpha_{2}-\delta}{\alpha_{2}-L}.
\end{equation}
\end{lem}

\noindent\textbf{Proof.} Fix a $t\in (0, T]$. If $\sqrt{t}<a_3$, let $n\ge 1$ satisfy
\begin{equation}
r^{n+1}a_1 = a_{2n+3} \le \sqrt{t} < a_{2n+1} = r^na_1.
\end{equation}
Then,
\begin{equation}
\f{\varrho a_{2n+1} - a_{2n+2}}{\sqrt{t}}\ge \varrho -p
\end{equation}
and
\begin{equation}
\f{a_{2n+2}}{\sqrt{t}}\le\frac{1}{q}.
\end{equation}
Therefore, by Lemma~\ref{lem:density},
\begin{equation}
\Prob{Y_t\in G} \ge 
\begin{cases}
\Prob{Y_{t}\in[a_{2n+2}, \varrho a_{2n+1}]}
\ge \Prob{\left|\mathcal{N}\right|\ge 1/q}\,\Prob{U\le \varrho-p}\!, & \text{if}\;\;\sqrt{t}<a_3, \\
\Prob{Y_{t}\in[a_2,\infty)}\ge \Prob{\left|\mathcal{N}\right|\ge a_2/\sqrt{t}}
\ge \Prob{\left|\mathcal{N}\right|\ge 1/q}\!, & \text{if}\;\;\sqrt{t}\ge a_3.
\end{cases}
\end{equation}

\smallskip

Since $c_3$ (appearing in the definition of $U$) stays bounded as $\alpha_2\downarrow1$, 
\begin{equation}
\liminf_{\alpha_2\downarrow 1}\,\Prob{Y_t\in G}
\ge \Prob{\left|\mathcal{N}\right|\ge 1/q}\,\liminf_{\alpha_2\downarrow 1} \Prob{U\le \varrho - p}
=:\iota > 0.
\end{equation}
Thus, for any $\alpha_2>1$ close enough to $1$,
\begin{equation}
\Prob{Y_t\in G}\ge\frac{\iota}{2}.
\end{equation}
Choosing
\begin{equation}
\delta = \alpha_2 -\f{\iota}{2}(\alpha_2- L)
\end{equation}
we obtain, for any $\alpha_2>1$ close enough to $1$, 
\begin{equation}
\Prob{Y_t\in G} \ge \f{\iota}{2} = \f{\alpha_2-\delta}{\alpha_2-L}\quad\text{and}\quad
\delta<1,
\end{equation}
and hence, \eqref{eq: final lemma}. \qed

\bigskip\bigskip\bigskip

%%%%%%%%%%%%%%%%%%%%%%%%%%%
\bibliographystyle{amsalpha}
\bibliography{biblio}
%%%%%%%%%%%%%%%%%%%%%%%%%%%

\bigskip\bigskip\bigskip

\end{document}